\documentclass{amsart}

\usepackage{palatino, epsfig, amsmath, amsthm, amssymb, mathpazo, mathrsfs}
\usepackage{eucal, eufrak}

\input diagrams.tex

\begin{document}

\newtheorem{theorem}{Theorem}[section]
\newtheorem{lemma}[theorem]{Lemma}
\newtheorem{corollary}[theorem]{Corollary}
\newtheorem{conjecture}[theorem]{Conjecture}
\newtheorem{question}[theorem]{Question}
\newtheorem{problem}[theorem]{Problem}
\newtheorem*{claim}{Claim}
\newtheorem*{criterion}{Criterion}
\newtheorem*{torus_thm}{Theorem A}

\theoremstyle{definition}
\newtheorem{definition}[theorem]{Definition}
\newtheorem{construction}[theorem]{Construction}
\newtheorem{notation}[theorem]{Notation}
\newtheorem{convention}[theorem]{Convention}

\theoremstyle{remark}
\newtheorem{remark}[theorem]{Remark}
\newtheorem{example}[theorem]{Example}

\def\A{\mathbb A}
\def\R{\mathbb R}
\def\H{\mathbb H}
\def\Z{\mathbb Z}
\def\C{\mathbb C}
\def\Q{\mathbb Q}
\def\B{\text{B}}
\def\E{\text{E}}
\def\F{\mathcal F}
\def\G{\mathcal G}
\def\P{\mathcal P}
\def\SL{\text{SL}}
\def\SO{\text{SO}}
\def\O{\text{O}}
\def\PSL{\text{PSL}}
\def\PGL{\text{PGL}}
\def\RP{\mathbb{RP}}
\def\CP{\mathbb{CP}}
\def\til{\widetilde}
\def\isom{\text{Isom}}
\def\homeo{\text{Homeo}}
\def\bhomeo{\text{BHomeo}}
\def\Gal{\text{Gal}}
\def\spin{\text{spin}}
\def\geo{\text{geo}}
\def\trace{\text{trace}}
\def\Ext{\text{Ext}}
\def\Hom{\text{Hom}}
\def\rot{\text{rot}}
\def\u{\text{univ}}

\title{Real places and torus bundles}
\author{Danny Calegari}
\address{Department of Mathematics \\ California Institute of Technology \\ Pasadena CA, 91125}

\date{2/16/2006, Version 0.15}

\begin{abstract}
If $M$ is a hyperbolic once-punctured torus bundle over $S^1$, then the trace field
of $M$ has no real places. 
\end{abstract}

\maketitle

\section{Introduction}

This paper studies a particular example of the interaction between topology and number theory
in the context of finite volume hyperbolic $3$-manifolds. If a {\em noncompact} hyperbolic
$3$-manifold $M$ is irreducible and
atoroidal, and homeomorphic to the interior of a compact $3$-manifold with 
torus boundary components, then $M$ admits a unique complete hyperbolic structure.
This structure determines a faithful representation from $\pi_1(M)$ into $\PSL(2,\C)$
which can be taken to have image in $\PSL(2,K)$ for some smallest number field $K$, 
called the {\em trace field} of $M$. See \cite{MacReid}, Theorem 4.2.3. 
(Note that the noncompactness of $M$ is important for the identification of $K$ with the
trace field of $\pi_1(M)$.)

It is an important question to study which number fields $K$ can arise from this
construction, and conversely, to understand the relationship between the topology
of $M$ and the algebra of $K$. This question seems to be wide open, and very
few nontrivial relationships are known.

In this paper we show by an elementary geometric argument that when $M$ is
a hyperbolic once-punctured torus bundle over $S^1$, the trace field $K$ has
no real places (i.e. it is not Galois conjugate to a subfield of $\R$).

\subsection{Statement of results}

In \S 2 we give an exposition of the relation between Euler classes, Stiefel-Whitney
classes, orientations, and boundary traces for hyperbolic surfaces and $3$-manifolds.
All this material is classical, but it seems that there is no explicit and thorough
account of it in the literature which takes our particular viewpoint. This discussion
takes up more than half the length of the paper --- we feel that its inclusion is
justified by its potential interest as a reference for the $3$-manifold community.

Note that some facts which might be otherwise somewhat obscure, are clarified
by a thorough discussion of this foundational material.
In particular, it is immediate from our viewpoint that for a hyperbolic
knot complement $M$, the trace of a longitude is $-2$, after lifting a geometric
representation $\pi_1(M) \to \PSL(2,\C)$ to $\SL(2,\C)$. 

\vskip 12pt

In \S 3 we specialize to once-punctured torus bundles. After a preliminary discussion
of trace fields and invariant trace fields, we prove our main theorem:

\begin{torus_thm}
Let $M_\phi$ be a hyperbolic oriented once-punctured torus bundle over $S^1$ with
monodromy $\phi$. Let $K$ denote the trace field of $M_\phi$, and $k$ the invariant trace field. 
Then:
\begin{enumerate}
\item{$K$ admits no real places}
\item{If $k$ admits a real place, then $H_1(M)$ contains 2-torsion}
\end{enumerate}
\end{torus_thm}

\vskip 12pt

In \S 4 we make some related observations. Most interesting is the observation
that if $M$ contains a {\em pseudo-rational surface} --- i.e. a surface $S$ with
all traces rational --- then such a surface has maximal Euler class with
respect to any $\PSL(2,\R)$ representation. In particular, a $3$-manifold containing
a pseudo-rational surface which is not Thurston norm minimizing in its homology class
(e.g. it might be separating) has a trace field with no real place. A nice corollary,
pointed out by A. Reid, is that this equality lets us construct examples of knot
complements which contain no totally geodesic immersed surfaces, for example, the knot
$8_{20}$ in Rolfsen's tables.

We also give an example showing that for every $n<0$ and every $m$ with $|m| \le n$ and
$m=n \mod 2$ there are incompressible surfaces $S_{n,m}$ with Euler characteristic
$n$ in some hyperbolic $3$-manifold $M$ with a real place $\sigma$ for which
$e_\sigma(S_{n,m}) = m$. The Milnor-Wood inequality implies that $|m| \le n$
is sharp, and the mod 2 condition is implied by the fact that the representations
are Galois conjugate into a geometric representation, so these are all the
possibilities which can arise. Somewhat amazingly, all these examples can be found in a 
single manifold $M$, the complement of the link $8_6^2$ in Rolfsen's tables.

\subsection{Acknowledgements}

Thanks to J. Button, N. Dunfield, O. Goodman, A. Reid and the anonymous referee
for comments, computations and corrections. 
Extra special thanks to A. Reid for nursing this paper through an
problematic gestation, and for the idea of Corollary~\ref{no_totally_geodesic}.

\section{Euler, Stiefel-Whitney, traces}\label{characteristic_class_section}

In what follows, we frequently deal with complete noncompact hyperbolic
surfaces $S$ with finite area. Such surfaces are homeomorphic to the interior
of a compact surface $\overline{S}$ with boundary. By abuse of notation, we will refer
to the boundary components of the (compactifying) surface $\overline{S}$ as the
boundary components of $S$.

For the convenience of the reader, we give a thorough 
exposition of the relationship between
$\SL,\PSL$, Stiefel-Whitney, Euler, and hyperbolic geometry. 
Related references are \cite{Goldman}, \cite{Culler}.

\subsection{Geometric representation}

Let $S$ be a noncompact orientable surface of finite type, and let 
$\phi:S \to S$ be a pseudo-Anosov homeomorphism. Then the mapping torus
$$M_\phi := S \times I /(s,1) \sim (\phi(s),0)$$
admits a complete hyperbolic structure. Corresponding to this hyperbolic
structure there is a discrete faithful representation
$$\rho_\geo:\pi_1(M_\phi) \to \PSL(2,\C)$$
which is unique up to conjugacy and orientation. We fix one such representation,
and call it the {\em geometric representation}.

\subsection{Quasifuchsian representations}
An orientable surface $S$ with negative Euler characteristic
itself admits a hyperbolic structure, and any such hyperbolic
structure determines some discrete faithful representation
$$\rho_S:\pi_1(S) \to \PSL(2,\R)$$
for which the image of each boundary curve is parabolic.

The space of discrete faithful representations from $\pi_1(S)$ into
$\PSL(2,\C)$ up to conjugacy and orientation, and for which the image of each
boundary curve maps to a parabolic element, is connected, and
contains an open dense subset homeomorphic to a ball consisting 
of {\em quasifuchsian} representations. For each quasifuchsian representation
$\rho_q:\pi_1(S) \to \PSL(2,\C)$ the ideal circle of the universal cover $\til{S}$
maps to a quasicircle in the sphere $\CP^1$. See \cite{Brock_Canary_Minsky} 
and \cite{MatTan} for background
and details on the theory of Kleinian groups.

\subsection{Nonorientable representations}

A surface $S$ might be orientable and yet admit representations into $\isom(\H^2)$
with nonorientable holonomy. The group of all isometries of $\H^2$ is $\PGL(2,\R)$ which
embeds into $\PGL(2,\C) \cong \PSL(2,\C)$.

\subsection{Classifying spaces}

We may embed both $\PSL(2,\R)$ and $\PGL(2,\R)$ into the group $\homeo(S^1)$ by considering
their action on the ideal circle $S^1_\infty$ of $\H^2$. The image of $\PSL(2,\R)$ is contained in
$\homeo^+(S^1)$ (the superscript $+$ denotes orientation preserving homeomorphisms).
These inclusion maps are homotopy equivalences.
The circle $S^1$ embeds in $S^2$ as the equator, and the group $\homeo(S^1)$ embeds in
$\homeo^+(S^2)$. We have two commutative diagrams
$$
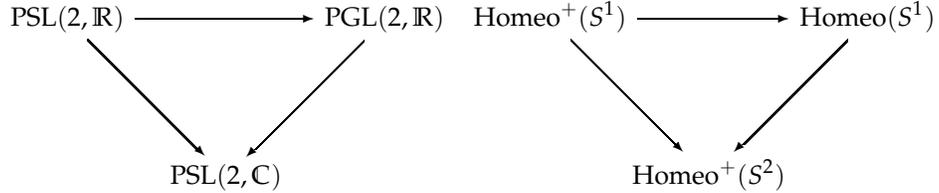
\begin{diagram}
\PSL(2,\R)  &  			& \rTo 				& 				&\PGL(2,\R) \\
						& \rdTo &							& \ldTo		&						\\
						&				& \PSL(2,\C)	&					&						\\
\end{diagram}
\; \; \; \;
\begin{diagram}
\homeo^+(S^1)  &  			& \rTo 				& 				&\homeo(S^1) \\
						& \rdTo &							& \ldTo		&						\\
						&				& \homeo^+(S^2)	&					&						\\
\end{diagram}$$
where the left diagram includes into the right diagram by a homotopy equivalence.

It is further true that the inclusions $$\homeo^+(S^1) \to \homeo^+(D^2)$$ and
$$\homeo^+(S^2) \to \homeo^+(D^3)$$ where $D^2,D^3$ are open balls,
obtained by coning to the center and throwing away the boundary, are homotopy equivalences.
The first case is straightforward; the second follows from Hatcher's proof of the
Smale conjecture \cite{Hatcher} and the (homotopy) equivalence of the categories
DIFF and TOP in dimension $3$.

For any group $G$, and any space $X$, a representation
$$\rho:\pi_1(X) \to G$$
induces a homotopy class of maps to the classifying space $\B G$.
There is a tautological $G$ bundle over $\B G$ called $\E G$ which pulls back
to a $G$-bundle over $X$ called $E_\rho$. Topologically, 
we form $E_\rho$ as the quotient bundle
$$E_\rho = \til{X} \times G / (x,g) \sim (\alpha(x),\rho(\alpha)(g))$$
where $\alpha$ ranges over elements of $\pi_1(X)$.

\subsection{Homotopy of $\B \homeo$}

The group $\homeo^+(S^1)$ is homotopy equivalent to the subgroup $\SO(2,\R) \approx S^1$ consisting of
rotations. $\homeo(S^1)$ is homotopy equivalent to the subgroup $\O(2,\R)$ which is homotopic to two
disjoint circles. $\homeo^+(S^2)$ is homotopy equivalent to $\SO(3,\R) \approx \RP^3$. It follows that we
can compute the homotopy groups of $\B\homeo$:
$$\pi_i(\B\homeo^+(S^1)) = \begin{cases}
\Z & \text{ if }i=2 \\
0 & \text{ otherwise } \\
\end{cases}
\; \; \;
\pi_i(\B\homeo(S^1)) = \begin{cases}
\Z/2\Z & \text{ if }i=1 \\
\Z & \text{ if }i=2 \\
0 & \text{ otherwise } \\
\end{cases}
$$
$$\pi_i(\B\homeo^+(S^2)) = \begin{cases}
\Z/2\Z & \text{ if }i=2 \\
0 & \text{ if }i=0,1,3 \\
\Z & \text{ if }i=4 \\
\end{cases}
$$
and is torsion for $i>4$.

Cohomology on $\B G$ pulls back to cohomology classes on $X$ which represent the first
obstruction to trivializing the bundle $E_\rho$. Since $\B\homeo^+(S^1)$ is a $K(\Z,2)$ and therefore
has the homotopy type of $\CP^\infty$, the cohomology ring is generated by a single
element $e \in H^2(\B\homeo^+(S^1);\Z)$. Identifying $\homeo^+(S^1)$
with $\PSL(2,\R)$ up to homotopy, this element represents the obstruction to lifting a
representation from $\PSL(2,\R)$ to $\til{\SL}(2,\R)$, its universal covering group.
The mod 2 reduction of $e$ is the obstruction to lifting to $\SL(2,\R)$.

Since $\B\homeo^+(S^2)$ is a $K(\Z/2\Z,2)$ below dimension $4$, only the class
$w \in H^2(\B\homeo^+(S^2);\Z/2\Z)$ is relevant to $3$-manifolds. Identifying $\homeo^+(S^2)$
with $\PSL(2,\C)$ up to homotopy, this represents the obstruction to lifting a
representation from $\PSL(2,\C)$ into $\SL(2,\C)$. If we include $\PSL(2,\R)$ into
$\PSL(2,\C)$ then we see that $w$ is the image of $e$ under $H^2(X;\Z) \to H^2(X;\Z/2\Z)$.

\subsection{$\B\homeo(S^1)$}

The short exact sequence 
$$0 \to \homeo^+(S^1) \to \homeo(S^1) \to \Z/2\Z \to 0$$
gives rise to a fibration of spaces
$$\B\homeo^+(S^1) \to \B\homeo(S^1) \to \RP^\infty$$
which exhibits $\B\homeo(S^1)$ up to homotopy as a twisted
$\CP^\infty$ bundle over $\RP^\infty$. The generator of $\pi_1(\B\homeo(S^1))$
acts on $\pi_2(\B\homeo(S^1))$ by multiplication by $-1$. Since this multiplication is trivial with
$\Z/2\Z$ coefficients, the $\Z/2\Z$ cohomology can be computed from the K\"unneth formula. In low
dimensions we get
$$H^i(\B\homeo(S^1);\Z/2\Z) =  \begin{cases}
\Z/2\Z & \text{ for } i=0,1 \\
\Z/2\Z \oplus \Z/2\Z & \text{ for } i=2 \\
\end{cases}$$
The generator in dimension $1$ is the orientation class $o$, and the generators in dimension $2$
are $o^2$ and $r$ which is the mod 2 reduction of the Euler class, and is obtained by pulling back
$w$ from $H^2(\B\homeo^+(S^2);\Z/2\Z)$.

\subsection{Relative bundles}

Suppose we have a pair of spaces $X,Y$ and a representation
$$\rho:\pi_1(X) \to G$$
such that $\rho|_{\pi_1(Y)}$ is trivial. Then the resulting map to the classifying space maps $Y$ to
the basepoint of $\B G$, and there is a well-defined relative homotopy class of pairs. This defines
a canonical trivialization of the restricted bundle $E_\rho|_Y$.
It follows that we can pull back reduced cohomology of $\B G$ to a relative class in $H^*(X,Y)$
which represents the obstruction to extending this trivialization of $E_\rho$ over $Y$ to all of $X$.

More generally, suppose $G$ contains a closed contractible subgroup $H$. Then the coset space $G/H$ is homotopy
equivalent to $G$, and the bundle $E_\rho$ can be replaced by a homotopy equivalent bundle $E_\rho/H$.
A representation $\rho:\pi_1(X) \to G$ for which $\rho(\pi_1(Y))$ is contained in $H$ defines a
canonical trivialization of $E_\rho/H$, and therefore a relative cohomology class in $H^*(X,Y)$.

In the groups $\PSL(2,\R)$ and $\PSL(2,\C)$ the maximal parabolic subgroups are contractible; these are
the $N$ subgroups with respect to the $KAN$ (or Iwasawa) decomposition of these groups.

In $\homeo^+(S^1)$ the stabilizer of a point is contractible. In $\homeo^+(S^2)$, the stabilizer of a point
is not contractible, being homotopic to $S^1$.

\subsection{Euler classes of surfaces, and relative Euler class}

For $S$ a closed orientable hyperbolic surface, the Euler class $e_{\rho_S}$ of the
representation 
$$\rho_S: \pi_1(S) \to \PSL(2,\R)$$ 
associated to any hyperbolic structure on $S$ satisfies
$$e_{\rho_S}([S]) = \pm \chi(S)$$
where $[S]$ represents the fundamental class in $H_2(S;\Z)$, and the sign depends on the
choice of orientation. Note in this
case that $E_{\rho_S}$ is isomorphic to the unit tangent bundle of $S$.

If $S$ is orientable and complete with finite area, but possibly with punctures,
then there is a relative fundamental class in $H_2(S,\partial S;\Z)$.
If $E_{\rho_S}$ represents the associated circle bundle, then $E_{\rho_S}$ restricts
to a finite union of tori fibering over the boundary components of $S$.

By our assumption, the holonomy around a boundary component is parabolic,
and has a (unique) fixed point in $S^1$. This fixed point suspends to a canonical
section of $E_{\rho_S}$ over $\partial S$, and defines a trivialization
of this restricted bundle. The {\em relative Euler class} of $\rho_S$ is the
element of $H^2(S,\partial S;\Z)$ which represents the obstruction to extending this
trivialization of $E_{\rho_S}|_{\partial S}$ over all of $E_{\rho_S}$.

\subsection{Geometric computation of Euler class}

In \cite{Thurston_circles}, Thurston defined a $2$-cocycle for a group $G$ acting
in an orientation-preserving way on $S^1$.
If $\sigma:G \to \homeo^+(S^1)$ is a representation,
choose a point $p \in S^1$ and for every triple $g_0,g_1,g_2 \in G$ define
$$c(g_0,g_1,g_2) =  \begin{cases}
1 & \text{if $g_0(p),g_1(p),g_2(p)$ are positively ordered} \\
-1 & \text{if $g_0(p),g_1(p),g_2(p)$ are negatively ordered} \\
0 & \text{if $g_0(p),g_1(p),g_2(p)$ are degenerate} \\
\end{cases}$$
Equivalently, think of $S^1$ as the ideal boundary of $\H^2$. Then for $g_0,g_1,g_2$ a triple
of elements in $G$, define $c(g_0,g_1,g_2)$ to be the (signed) area of the ideal triangle
with vertices at the $g_i(p)$, divided by $\pi$.

From this geometric definition, it is easy to see that $c$ is coclosed, and that the
cohomology class it defines is independent of $p$. 
Suppose that
$G = \pi_1(M)$ and $H_i < G$ is $\pi_1(\partial M_i)$ for each component 
$\partial M_i$ of $\partial M$. If each $H_i$ has a fixed point in $S^1$, then 
for some choice of basepoint, $c$ is identically zero on $H_i$. 
It follows that $c$ defines a relative class
in $H^2(M,\partial M;\Z)$.

\begin{lemma}\label{Thurston_and_Euler}
The Thurston cocycle $c$ is related to the Euler class $e$ by
$$[c] = 2[e] \text{ in } H^2(M,\partial M;\Z)$$
\end{lemma}

For a proof, see e.g. \cite{Jekel}. For more details, and related constructions,
see \cite{Ghys} or \cite{Calegari_euler}.

\subsection{Stiefel-Whitney class and $\SL(2,\C)$}

For a closed hyperbolic $3$-manifold, the $S^2$ bundle coming from the
geometric representation $\rho_\geo$ is isomorphic to the unit tangent bundle (to see
this, use the exponential map). The
pullback of the generator of $H^2(\bhomeo^+(S^2);\Z/2\Z)$ can therefore be identified with
the second Stiefel-Whitney class of $M$. If $M$ is orientable, $TM$ is parallelizable,
so this class must vanish. 

If $M$ is noncompact with finite volume,
the geometric representation on each boundary torus group is parabolic, with a
fixed point $p \in \CP^1$. The subgroup $N_p < \PSL(2,\C)$ of parabolic elements
fixing $p$ is contractible, so there is a trivialization of
the associated bundle over $\partial M$ and we get a relative second Stiefel-Whitney
class
$$w_\geo \in H^2(M,\partial M;\Z/2\Z)$$
which might not be trivial.

Note that the image of $w_\geo$ in $H^2(M;\Z/2\Z)$ is the ordinary second Stiefel-Whitney class
which is always trivial for orientable $M$, as above. 
It follows from the long exact sequence in cohomology that 
$w_\geo$ is the image of a distinguished
class in $H^1(\partial M;\Z/2\Z)/i^*H^1(M;\Z/2\Z)$, where $i^*$ denotes
the homomorphism induced by the inclusion map $i:\partial M \to M$.

Since the ordinary second Stiefel-Whitney class vanishes,
the geometric representation lifts to the double cover 
$$\hat{\rho}_\geo:\pi_1(M) \to \SL(2,\C)$$
Since every element of $\pi_1(\partial M)$ is parabolic, the preimages either have trace
$2$ or $-2$. Different choices of lift to $\SL(2,\C)$ might change the value on elements of
$H_1(\partial M;\Z/2\Z)$ which are in the image of $H_1(M;\Z/2\Z)$ but do not change the
values on elements which are homologically nontrivial in $\partial M$ but bound in $M$ (mod 2).
In particular, if $\alpha \in \pi_1(\partial M)$ represents
zero in $H_1(M;\Z/2\Z)$, then the trace is {\em independent} of the choice of lift.
We deduce that if $\trace(\hat{\rho}_\geo(\alpha)) = -2$ where $\alpha = \partial S$
for some properly embedded $S \subset M$ then $w_\geo$ is nontrivial.

In general, for any representation
$$\rho:\pi_1(M) \to \PSL(2,\C)$$
which sends $\pi_1(\partial M)$ to parabolic elements, we get a relative class 
$$w_\rho \in H^2(M,\partial M;\Z/2\Z)$$
If $\rho$ lifts to $\hat{\rho}:\pi_1(M) \to \SL(2,\C)$ then
the image of $w_\rho$ in $H^2(M;\Z/2\Z)$ is zero.
If $\trace(\hat{\rho}(\alpha)) = -2$ for some $\alpha \in \pi_1(\partial M)$
which represents a trivial class in $H_1(M;\Z/2\Z)$
then $w_\rho$ is nontrivial. We summarize this as a lemma:

\begin{lemma}\label{Stiefel_nontrivial}
Suppose that $M$ is a compact orientable $3$-manifold with torus boundary components,
and suppose $\rho$ is some representation
$$\rho:\pi_1(M) \to \PSL(2,\C)$$
for which elements of $\pi_1(\partial M)$ map to parabolic transformations. Suppose further
that $\rho$ lifts to 
$$\hat{\rho}:\pi_1(M) \to \SL(2,\C)$$
Let $w_\rho \in H^2(M,\partial M;\Z/2\Z)$ denote the relative Stiefel-Whitney
class of $\rho$, where the trivialization on $\partial M$ comes from the
fixed point of the corresponding parabolic subgroup.
Further, suppose there is some $\alpha \in \pi_1(\partial M)$
representing zero in $H_1(M;\Z/2\Z)$ for which 
$$\trace(\hat{\rho}(\alpha)) = -2$$
Then $w_\rho$ is nontrivial.
\end{lemma}

\subsection{$\chi$ and $w$}

Let $M$ be a manifold containing an incompressible surface $S$ with a single boundary component.
It turns out that there is a very simple relationship between $\chi(S)$ and $w_\geo([S])$; in fact,
$w_\geo([S])$ is just the mod 2 reduction of $\chi(S)$. 

The simplest way to see this is topological. The value of $w_\geo([S])$ depends only on the 
topology of the $S^2$ bundle over $S$ obtained from the representation 
$$\rho_\geo|_S:\pi_1(S) \to \PSL(2,\C) < \homeo^+(S^2)$$
By using the exponential map, 
we may identify this bundle with the restriction of the unit tangent bundle
$UTM|_S$. In particular, the value of $w_\geo([S])$ only depends on the topology of a tubular
neighborhood of $S$ in $M$. If $S$ is two-sided and embedded, then this tubular neighborhood
is just a product $S \times I$ and therefore the value of $w_\geo([S])$ depends only on $S$.
Now, for a geometric representation $\rho_S|_S$ coming from a hyperbolic structure on $S$,
we have already seen that $w_S([S]) =\chi(S) \text{ mod } 2$.
In particular, we have the following lemma:

\begin{lemma}\label{Stiefel_class}
Let $M$ be a complete finite-volume orientable hyperbolic $3$-manifold.
If $S \subset M$ is incompressible and orientable,
then
$$w_\geo([S]) = \chi(S) \mod 2$$
where $[S]$ represents the fundamental class in $H_2(S,\partial S;\Z/2\Z)$.
\end{lemma}

The following corollary is folklore, and seems to have been observed first by W. Thurston, at
least for Seifert surfaces of knots in $S^3$.
There appears to be some confusion in the literature about whether it
is well-known, and therefore we state it for completeness:

\begin{corollary}
Let $M$ be a complete noncompact orientable hyperbolic $3$-manifold, and
let $\alpha \in \partial M$ be the boundary of a $2$-sided incompressible surface
$S \subset M$. If $\hat{\rho}_\geo:\pi_1(M) \to \SL(2,\C)$ is any lift of the
geometric representation, then
$$\trace(\hat{\rho}_\geo(\alpha)) = -2$$
\end{corollary}

It follows that the longitude of any knot has trace $-2$. This answers Question 6.2
in \cite{Cooper_Long}.

\vskip 12pt

We remark that by the proof of the Ending Lamination Conjecture \cite{Brock_Canary_Minsky},
for $S$ incompressible, the geometric representation $\rho_\geo|_S$ coming from the 
hyperbolic structure on $M$ and a geometric representation $\rho_S|_S$ coming from a hyperbolic
structure on $S$ are homotopic through representations of $\pi_1(S)$ into $\PSL(2,\C)$ 
which send boundary curves to parabolic elements, 
since $\rho_\geo|_S$ is always in the closure of the space of quasifuchsian representations.
This gives another, more high-powered proof of Lemma~\ref{Stiefel_class}.

\section{Torus bundles}\label{torus_section}

\subsection{Number fields}

If $M$ is a complete hyperbolic $3$-manifold of finite volume, we have
the corresponding geometric representation
$$\rho_\geo:\pi_1(M) \to \PSL(2,\C)$$
The {\em trace field} $K$ of $M$ is the field generated by the traces of
$\rho_\geo(\alpha)$, as $\alpha$ varies over the elements of $\pi_1(M)$. Of course,
the trace of an element in $\PSL(2,\C)$ is only determined up to sign, but the
field generated by these elements is independent of sign.

It turns out that the trace field $K$ is always a {\em number field}
(i.e. some finite algebraic extension of $\Q$) and as remarked in the
introduction, {\em for $M$ noncompact}, $\rho_\geo$ can be conjugated into $\PSL(2,K)$. 

The {\em invariant trace field} is the subfield $k$ of $K$ generated by the
squares of the traces of $\rho_\geo(\alpha)$ as above. In general, $k$ and $K$ are not equal,
and the degree satisfies $[K : k] = 2^n$ for some $n$. The field $k$ is an invariant
of the commensurability class of $M$, where two hyperbolic manifolds $M,N$ are said to
be commensurable if they have a common finite cover.

See \cite{MacReid} for details and proofs.

\subsection{Action of $\Gal(L/\Q)$}

Let $L$ denote the Galois closure of $K$ in $\C$. If $p(x)$ is the minimal polynomial
of a generating element of $K$, then $L$ is obtained from $\Q$ by adjoining all roots
of $p(x)$.

The Galois group $\Gal(L/\Q)$ of $L$ over $\Q$ acts on $L$ by field
automorphisms, conjugating $K$ into different subfields of $\C$.

The various (Galois conjugate) embeddings of $K$ into $\C$ are called
{\em places}, and can be real or complex. Complex places come in pairs,
interchanged by complex conjugation. If the degree $[K:\Q]=d$, then
$$d = r_1 + 2r_2$$
where $r_1$ is the number of real places, and $r_2$ is the number of conjugate
pairs of complex places.

In general, if
$$\rho:\pi_1(M) \to \PSL(2,K)$$
is some representation which is parabolic on $\partial M$, and if the relative
Stiefel-Whitney class $w_\rho$ has zero image in $H^2(M;\Z/2\Z)$, then
$\rho$ lifts to $\hat{\rho}$, and parabolic elements lift to elements of $\SL(2,K)$ with
trace equal to $\pm 2$. If $\rho_\sigma$ is obtained from $\rho$ by Galois conjugating $K$
to $K_\sigma$, then $\rho_\sigma$ lifts to $\SL(2,K_\sigma)$, so the relative Stiefel-Whitney
class $w_{\rho_\sigma}$ has zero image in $H^2(M;\Z/2\Z)$. Moreover, since $\pm 2$ are in the
fixed field of $\sigma$, we have equality
$$\hat{\rho}_\sigma(\alpha) = \hat{\rho}(\alpha)$$
Since the relative class is determined by these traces,
we have equality
$$w_\rho = w_{\rho_\sigma}$$

\begin{remark}
The invariance of $w_\rho$ under the Galois group implies that $w_\rho$ can be pulled back from the
cohomology of $\B\PSL(2,\overline{\Q})$ where $\overline{\Q}$ has the discrete topology.
\end{remark}

\subsection{Real places for $K$ and $k$}

A real place for $K$ determines an embedding of $\pi_1(M)$ in $\PSL(2,\R)$, which lifts to
$\SL(2,\R)$, by the vanishing of the second Stiefel-Whitney class for a $3$-manifold, and
the invariance of $w$ under the action of the group $\Gal(L/\Q)$.

Let $\Gamma$ be the group generated by squares of elements of $\pi_1(M)$ where $M$ is noncompact
as before. A real place for $k$ determines an embedding of $\Gamma$ into $\PSL(2,\R)$
which extends to an embedding of $\pi_1(M)$ into $\PSL(2,\C)$ for which every
element of $\pi_1(M)$ has a trace which is real or pure imaginary.

Since the square of every element $\alpha$ of $\pi_1(M)$ stabilizes $\H^2$ in $\H^3$, it follows
that $\alpha$ either fixes $\H^2$ (possibly reversing orientation) or takes it to an orthogonal
copy of $\H^2$. In the second case, the square $\alpha^2$ takes $\H^2$ to itself by an
orientation-reversing isometry, contrary to the fact that $\alpha^2$ is in $\PSL(2,\R)$ by hypothesis.
It follows that $\pi_1(M)$ preserves $\H^2$, and therefore has image in $\PGL(2,\R)$.

Again, this discussion depends on the noncompactness of $M$. In general, for noncompact $M$,
a similar argument shows that 
the representation of $\pi_1(M)$ can always be conjugated into $\PGL(2,k)$ where $k$ is
the invariant trace field. This fact is implicit e.g. in \cite{Neumann_Reid}, page 278.

\subsection{Homology of torus bundles}

Let $M_\phi$ be a hyperbolic once-punctured torus bundle with monodromy $\phi$.
Then $\phi$ induces an automorphism on $H_1(T;\Z)$, represented by some matrix
$$\phi \sim \begin{pmatrix}
a & b\\
c & d\\
\end{pmatrix}$$
Since $\phi$ is pseudo-Anosov, the trace $a+d$ satisfies $|a+d|>2$.
Then $H_1(M_\phi)$ is isomorphic to the kernel of the map
$$\begin{pmatrix}
a-1 & b & 0 \\
c & d-1 & 0 \\
0 & 0 & 0 \\
\end{pmatrix} \;: \; \Z^3 \to \Z^3$$
which has rank $1$, and therefore $H_2(M_\phi,\partial M_\phi;\Z)$ is isomorphic to $\Z$,
generated by the relative class of the fiber.

\subsection{Torus bundles and $2$-torsion}

Let $M_\phi$ be a hyperbolic surface bundle over $S^1$ with fiber a once-punctured
torus $T$. The following theorem relates topology, homological algebra, and
number theory:

\begin{torus_thm}
Let $M_\phi$ be a hyperbolic oriented once-punctured torus bundle over $S^1$ with
monodromy $\phi$. Let $K$ denote the trace field of $M_\phi$, and $k$ the invariant trace field. 
Then:
\begin{enumerate}
\item{$K$ admits no real places}
\item{If $k$ admits a real place, then $H_1(M)$ contains 2-torsion}
\end{enumerate}
\end{torus_thm}
\begin{proof}
We suppose after conjugating that the image of $\rho_\geo$ lies in $\PSL(2,K)$.
We know that $\rho_\geo$ lifts to 
$$\hat{\rho}_\geo:\pi_1(M_\phi) \to \SL(2,K)$$
Moreover, if $T$ denotes the fiber of the fibration, and $A,B$ are standard (free)
generators for $\pi_1(T)$, then
$$\text{trace }\hat{\rho}_\geo([A,B]) = -2$$
as in Lemma~\ref{Stiefel_class}.

Suppose $\sigma:K \to \R$ is a real place, and let 
$$\hat{\rho}_\sigma:\pi_1(M_\phi) \to \SL(2,\R)$$
be obtained by Galois conjugating $K$ into $\R$.
Let $\rho_\sigma:\pi_1(M_\phi) \to \PSL(2,\R)$ be obtained by composing $\hat{\rho}_\sigma$
with the covering map $\SL(2,\R) \to \PSL(2,\R)$.

We denote the relative Euler class of $\rho_\sigma$ by $e$:
$$e \in H^2(M_\phi,\partial M_\phi;\Z)$$
Since $w_\geo = w_\sigma$ we must have that $e([T])$ is odd.

Triangulate $T$ by two ideal triangles $\Delta_1,\Delta_2$. The representation
$\rho_\sigma$ determines a developing map from the universal cover $\til{T}$ to
$\H^2$
$$d:\til{T} \to \H^2$$ 
If the image of both triangles has the same orientation, then the developing
map is a homeomorphism, and we obtain a complete hyperbolic structure on $T$ which is
invariant under $\phi$. But this implies that $\phi$ has finite order, which is
incompatible with the existence of a complete hyperbolic structure on $M_\phi$.
It follows that the orientations on the images of the $\Delta_i$ disagree.

By Thurston's formula for $2e$ (Lemma~\ref{Thurston_and_Euler}), we have
$$2e([T]) = \sum_i \text{sign of orientation on } d(\Delta_i) = 0$$
This gives a contradiction, and shows that $K$ has no real place.

\vskip 12pt

Now if $k$ admits a real place, then there is some $\sigma: K \to \C$ such that
$$\rho_\sigma:\pi_1(M_\phi) \to \PGL(2,\R)$$
As above, we get a developing map from $\til{T}$ to $\H^2$ for which the orientations on
the ideal triangles must disagree, and the rational relative Euler class of the action must vanish.
Since $\rho_\sigma$ is conjugate into $\PGL(2,\R)$ but not $\PSL(2,\R)$ the orientation class
$o_\sigma \in H^1(M_\phi;\Z/2\Z)$ must be nontrivial. In fact, since 
the traces of elements of $\pi_1(\partial M_\phi)$
are $\pm 2$, boundary elements map to the subgroup $\PSL(2,\R)$, and therefore the orientation
class $o_\sigma$ is a nontrivial class in $H^1(M_\phi,\partial M_\phi;\Z/2\Z)$.

Since $H_1(M,\partial M;\Z)$ is torsion for a punctured torus bundle, we are done.
\end{proof}

\begin{remark}
Note that any field of odd degree admits a real place; in particular, the degree of $K$ is always even.
\end{remark}

\begin{remark}
If $M$ is a (compact or noncompact) hyperbolic surface bundle, and $S$ is any fiber, then
the same argument shows that if $K$ has a real place $\sigma:K \to \R$, then
$$|e_\sigma(S)| < -\chi(S), \; e_\sigma(S) = \chi(S) \mod 2$$
\end{remark}

\begin{remark}
J. Button has studied trace fields of punctured torus bundles with monodromy of the
form $L^{-1}R^{-n}$ in \cite{Button}. He showed for positive $n \equiv 2 \mod 4$
and for all odd $n$ that the invariant trace field of $\pi_1(M_n)$ has no real places.
For $n$ odd, the homology of $M_n$ has no $2$-torsion, but for $n$ even and not
divisible by $4$, this does not follow from Theorem A, but rather from an
explicit computation.

One might further ask whether every punctured torus bundle whose invariant trace
field has a real place has $4$-torsion in $H_1$. In fact, we posed exactly this
question in an earlier version of this paper. J. Button has found a counterexample
to this question: the census manifold s299 is a once-punctured torus bundle with monodromy
$-R^4L^2$ whose invariant trace field $k$ has degree 3, and whose trace field $K$ has degree 12,
and which has first homology $\Z \oplus \Z/2\Z \oplus \Z/6\Z$.
\end{remark}

\begin{remark}
In \cite{Goldman}, W. Goldman characterizes geometric representations
of once punctured torus groups amongst all $\PSL(2,\R)$ representations
in terms of trace data. This gives an alternate proof of
the first part of Theorem A, without using the geometric formula for the Euler class.
\end{remark}

\begin{remark}
Part (2) of Theorem A also follows from Corollary 2.3 of \cite{Neumann_Reid}.
\end{remark}

\begin{example}
Amongst the cusped manifolds in the Hodgson--Weeks census
(see \cite{Weeks}), 
m039 is a torus bundle with monodromy $RL^4$,
$H_1 = \Z \oplus \Z/4\Z$ and invariant trace field with minimal polynomial
$x^3 - x^2 + x + 1$. It has a degree $2$ cover v3225 for which this is the
trace field; this cover is fibered with fiber a twice-punctured torus, so
necessarily the Euler class of the representation associated to the real
place must vanish on this fiber.

Some other examples:
m040 is a torus bundle with monodromy $-RL^4$,
$H_1 = \Z \oplus \Z/8\Z$ and invariant trace field with minimal polynomial
$x^3 - x^2 + x + 1$,
and v2231 is a torus bundle with monodromy $RL^2RL^3$,
$H_1 = \Z \oplus \Z/16\Z$ and invariant trace field with minimal polynomial
$x^7 - 3x^5 - 2x^3 - 2x^2 + 4x - 2$.
The invariant trace fields were found with the help of the program {\tt snap}
(\cite{snap}).
\end{example}

\section{Inequalities for the Euler class}

\subsection{Thurston norm}

We have seen from \S~\ref{characteristic_class_section} and \S~\ref{torus_section}
that 
$$|e_\sigma(S)| < -\chi(S)$$
for $S$ a fiber of $M$, and
$$e_\sigma(S) = \chi(S) \mod 2$$
for any incompressible surface $S$,
whenever $\sigma:K \to \R$ is a real place.

In \cite{Thurston_norm}, Thurston introduced a norm on $H_2(M,\partial M;\R)$ for $M$
irreducible and atoroidal. For a homology class $[S]$, the norm satisfies
$$\|[S]\| = \inf_S -\chi(S)$$
where the infimum is taken over all (possibly disconnected) representatives $S$ of
$[S]$ with no spherical components.

A generalization of this norm, due to Gromov, measures a similar complexity amongst all
{\em immersed} surfaces with no spherical components representing a given homology class.
A theorem of Gabai (\cite{Gabai})
shows that these two norms are equal (after a suitable normalization);
i.e. any immersed surface may be replaced by an embedded surface of no larger norm.

The key properties of the norm $\|\cdot \|$ are that the unit ball $\P(M)$ is a finite sided
polyhedron, whose vertices are {\em rational}, and that there are a finite (possibly empty)
collection of top dimensional faces $Q_i$ with the property that the integral
homology classes $[S]$ representing
fibrations of $M$ over $S^1$ are exactly those whose (positive) projective rays intersect
the interiors of the $Q_i$. Such $Q_i$ are called {\em fibered faces} of $\P(M)$.

Our estimate implies the following:
\begin{theorem}\label{norm_ball}
Let $M$ be a cusped hyperbolic $3$-manifold, and 
suppose $\sigma:K \to \R$ is a real place with associated
relative Euler class $e_\sigma$. Then for every fibered face of $\P(M)$ there
is a vertex $V_i$ such that $e_\sigma(V_i) \ne \|V_i\|$. Similarly, there is
a vertex $V_j$ such that $e_\sigma(V_j) \ne -\|V_j\|$.
\end{theorem}

\subsection{Pseudo-rational surfaces}

In \cite{LongReid}, Long and Reid define a {\em pseudomodular} surface to
be one whose cusp set is contained in $\Q$.

We alter their definition slightly to adapt it to our context:

\begin{definition}
A subgroup $\Gamma < \PSL(2,\R)$ is {\em pseudo-rational} if the traces of all elements
are contained in $\Q$.
\end{definition}

A discrete finite covolume pseudo-rational subgroup acts on $\H^2$ with quotient a
pseudo-rational surface $S$. Such a surface in a hyperbolic $3$-manifold is necessarily
totally geodesic. Note that for us, pseudo-rational surfaces are always orientable.

\begin{example}
A thrice-punctured sphere is a (pseudo)-rational surface.
\end{example}

If $K$ is the trace field of $M$ and $\sigma:K \to \R$ is a real place, the traces
of a pseudo-rational subsurface do not change. It follows that 
$$e_\sigma([S]) = \pm \chi(S)$$
for any pseudo-rational surface, and any real place $\sigma$. Consequently, we have the
following corollary:

\begin{theorem}\label{pseudo_theorem}
Let $M$ be a cusped hyperbolic $3$-manifold, 
and suppose $S \subset M$ is a pseudo-rational
surface (possibly immersed). 
If $S$ is not (Gromov or Thurston) norm minimizing in its homology class, $K$ has no real places. 
\end{theorem}

In particular, if $M$ contains a separating pseudo-rational surface, its trace field
has no real places.

\begin{example}
One method of constructing (pseudo)-rational surfaces is by covering thrice punctured
spheres. A thrice-punctured sphere is always homologically essential, and therefore
so is its preimage in a finite cover. But if such a finite cover has suitable symmetries,
one might be able to find a (low genus)
surface with a $2$-fold orientation-reversing fixed-point free
symmetry, in the same homology class as the pseudo-rational surface. 
One can then cut along such a surface, and reglue
the resulting boundary components to themselves to get a new manifold with the
same trace field as the old, in which the pseudo-rational surface is homologically
trivial.
\end{example}

A nice application of Theorem~\ref{pseudo_theorem} is the following Corollary,
which was suggested by A. Reid:

\begin{corollary}\label{no_totally_geodesic}
Let $M$ be a fibered knot complement in a rational homology sphere
whose trace field $K$ has odd prime degree. Then $M$ does not contain an immersed totally
geodesic surface.
\end{corollary}
\begin{proof}
Since $K$ has prime degree, it has no proper subfields other than $\Q$; in particular, any 
immersed totally geodesic surface $S$ has rational traces.
Since $M$ is a knot complement in a rational homology sphere, its rational second homology is
$1$ dimensional. Since it is fibered, the rational homology is generated by the fiber $F$.
It follows that $[S] = n[F]$ in homology for some nonzero integer $n$. Since $F$ is a fiber of
a fibration, it is Thurston (and Gromov) norm-minimizing, and therefore
$$-\chi(S) \ge -|n|\chi(F)$$
Let $\sigma:K \to \R$ be a real place, and let $e_\sigma$ be the associated relative Euler class.
Then we have
$$|e_\sigma(S)| = |n|\cdot |e_\sigma(F)| < -|n|\chi(F) \le -\chi(S)$$
contrary to Theorem~\ref{pseudo_theorem}.
\end{proof}

For example, the knot $8_{20}$ in \cite{Rolfsen} (the complement is m222 in the census) 
is fibered, and has trace field generated by a root of
$x^5 - x^4 + x^3 + 2x^2 - 2x + 1$ which has degree $5$.

\subsection{Realizing Euler classes}

If $\Sigma$ is a closed, orientable surface of genus $g\ge 2$, Goldman
\cite{Goldman} showed that the $\PSL(2,\R)$
representation variety of $\pi_1(\Sigma)$ has $4g-3$ components, indexed by
values of the Euler class on $\Sigma$ satisfying
$$|e([\Sigma])| \le -\chi(\Sigma)$$

The Milnor--Wood inequality (c.f. \cite{Milnor}, \cite{Wood})
says that one cannot do better, even amongst
$\homeo^+(S^1)$ representations:
\begin{theorem}[Milnor--Wood]
Let $\Sigma$ be a closed surface of genus at least $1$, and let
$\rho:\pi_1(\Sigma) \to \homeo^+(S^1)$ be a representation with Euler
class $e_\rho$. Then
$$|e_\rho([\Sigma])| \le -\chi(\Sigma)$$
\end{theorem}

Similar theorems hold for surfaces with boundary, where one considers relative Euler classes.

For representations $\rho_\sigma$ coming from real places $\sigma$ of trace fields $K$,
we have the additional constraint that $e_\sigma([\Sigma]) = \chi(\Sigma) \mod 2$.
Modulo this constraint, we will see how to construct simple examples which realize
every possible compatible combination of Euler characteristic and Euler class.

\begin{definition}
Let $M$ be a manifold, and $\P(M)$ the unit ball of the Thurston norm. A
{\em big diamond} is a symmetrical $4$-gon $D \subset \P(M)$ which is the
intersection of $\P(M)$ with a two-dimensional plane $\pi$, and whose vertices
are integer lattice points which generate the lattice of integral points in $\pi$.
\end{definition}

Since the norm of every integer lattice point is at least $1$, a ``big diamond" is
as big as possible, hence the name. Notice too that only cusped manifolds can have
big diamonds in $\P(M)$.

\begin{theorem}
Let $M$ be a cusped hyperbolic $3$-manifold. Suppose the trace field $K$ has a real
place $\sigma$ with associated relative Euler class $e_\sigma$, 
and suppose further that the unit ball in the Thurston norm $\P(M)$ contains a big diamond
$D$. Then for every integer $n<0$ and every integer $m$ with $|m|\le -n$ and $n=m \mod 2$
there is an immersed incompressible connected surface  $S_{n,m}$ in $M$ satisfying
$$\chi(S_{n,m}) = n, \; e_\sigma([S_{n,m}]) = m$$
\end{theorem}
\begin{proof}
Let $V_1,V_2$ be surfaces representing the vertices of the big diamond $D$.
Then $\chi(V_1)=\chi(V_2) = -1$ and therefore $|e_\sigma(V_1)| = |e_\sigma(V_2)| = 1$.
After replacing $V_1$ and/or $V_2$ with their negatives if necessary, we can
assume 
$$e_\sigma(V_1) = 1, \; e_\sigma(V_2) = -1$$
For $p,q \ge 1$ let $V_{p,q}$ denote the Thurston norm-minimizing surface representing
the homology class $p[V_1] + q[V_2]$. Since $D$ is a diamond, we have
$$\chi(V_{p,q}) = -p-q$$
Since $e_\sigma$ is linear, we have
$$e_\sigma(V_{p,q}) = p-q$$
If $p,q$ are coprime, then $V_{p,q}$ is represented by a connected surface. Otherwise,
we have $p = ap',q=aq'$ for some $a>1$ where $p',q'$ are coprime.
Then $\pi_1(V_{p',q'})$ has a subgroup of index $a$ which gives a connected
incompressible immersed surface in $M$ with Euler characteristic $-p-q$ and
Euler class $p-q$. Together with finite index subgroups of $\pi_1(V_1),\pi_1(V_2)$,
this shows that every possibility is realized.
\end{proof}

\begin{example}
The link $8_6^2$ in Rolfsen's tables \cite{Rolfsen} has a complement whose unit ball
in the Thurston norm is a big diamond, 
and has trace field generated by a root of $x^3 -x^2 + 3x - 2$, which has a real
place because the degree is odd (thanks to N. Dunfield for finding this example.)
\end{example}


\begin{thebibliography}{99}
\bibitem{Brock_Canary_Minsky}
	J. Brock, R. Canary and Y. Minsky,
	\emph{The classification of Kleinian surface groups, II: The Ending Lamination Conjecture},
	eprint, math.GT/0412006
\bibitem{Button}
	J. Button,
	\emph{Invariant trace fields of once-punctured torus bundles},
	Kodai Math. J. {\bf 28} (2005), no. 1, 181--190
\bibitem{Calegari_euler}
	D. Calegari,
  \emph{Circular groups, planar groups and the Euler class},
  Geom. Topol. Monog. {\bf 7} (Proceedings of the Casson Fest) (2004), 431--491
\bibitem{Cooper_Long}
	D. Cooper and D. Long,
	\emph{Remarks on the $A$-polynomial of a knot},
	J. Knot Theory Ramifications {\bf 5} (1996), no. 5, 609--628
\bibitem{Culler}
	M. Culler,
	\emph{Lifting representations to covering groups},
	Adv. in Math. {\bf 59} (1986), no. 1, 64--70
\bibitem{Gabai}
	D. Gabai,
	\emph{Foliations and the topology of $3$-manifolds},
	J. Diff. Geom. {\bf 18} (1983), no. 3, 445--503
\bibitem{Ghys}
	E. Ghys,
	\emph{Groupes d'hom\'eomorphismes du cercle et cohomologie born\'ee},
	{\it The Lefschetz centennial conference, Part III (Mexico City, 1984)}, 81--106
	Contemp. Math., 58, III Amer. Math. Soc., Providence, RI, 1987
\bibitem{Goldman}
	W. Goldman,
	\emph{Topological components of spaces of representations},
	Invent. Math. {\bf 93} (1988), no. 3, 557--607
\bibitem{snap}
	O. Goodman, Snap, available at {\tt http://www.ms.unimelb.edu.au/\verb+~snap+/}
\bibitem{Hatcher}
	A. Hatcher,
	\emph{A proof of a Smale conjecture, ${\rm Diff}(S^3) \simeq {\rm O}(4)$},
	Ann. of Math. (2) {\bf 117} (1983), no. 3, 553--607
\bibitem{Jekel}
	S. Jekel,
	\emph{A simplicial formula and bound for the Euler class},
	Israel J. Math. {\bf 66} (1989), no. 1-3, 247--259
\bibitem{LongReid}
	D. Long and A. Reid,
	\emph{Pseudomodular surfaces},
	J. Reine Angew. Math. {\bf 552} (2002), 77--100
\bibitem{MacReid}
	C. Maclachlan and A. Reid,
	\emph{The arithmetic of hyperbolic 3-manifolds},
	Graduate Texts in Mathematics 219
	Springer-Verlag, New York, 2003
\bibitem{MatTan}
	K. Matsuzaki and M. Taniguchi,
	\emph{Hyperbolic manifolds and Kleinian groups},
	Oxford Mathematical Monographs, Oxford Science Publications
	The Clarendon Press, Oxford University Press, New York, 1998
\bibitem{Milnor}
	J. Milnor,
	\emph{On the existence of a connection with curvature zero},
	Comment. Math. Helv. {\bf 32} (1958), 215--223
\bibitem{Neumann_Reid}
	W. Neumann and A. Reid,
	\emph{Arithmetic of hyperbolic manifolds},
	Topology '90 (Columbus, OH, 1990), 273--310,
	Ohio State Univ. Math. Res. Inst. Publ., 1,
	de Gruyter, Berlin, 1992.
\bibitem{Rolfsen}
	D. Rolfsen,
	\emph{Knots and Links},
	Corrected reprint of the 1976 original Mathematics Lecture Series, 7. 
	Publish or Perish, Inc., Houston, TX, 1990
\bibitem{Thurston_norm}
	W. Thurston,
	\emph{A norm for the homology of $3$-manifolds},
	Mem. Amer. Math. Soc. {\bf 59} (1986), no. 339, i--vi and 99--130
\bibitem{Thurston_circles}
	W. Thurston,
	\emph{Three-manifolds, foliations and circles II},
	privately circulated preprint, (1997)
\bibitem{Weeks}
	J. Weeks, SnapPea, available at {\tt http://www.geometrygames.org/SnapPea/}
\bibitem{Wood}
	J. Wood,
	\emph{Bundles with totally disconnected structure group},
	Comment. Math. Helv. {\bf 46} (1971), 257--273
\end{thebibliography}
\end{document}